\documentclass[10pt]{amsart}
\usepackage{amssymb,amsmath,graphicx}

\theoremstyle{definition}
\newtheorem{theorem}{Theorem}[section]
\newtheorem{definition}{Definition}[section]
\newtheorem{proposition}{Proposition}[section]
\newtheorem{lemma}{Lemma}[section]

\newtheorem{corollary}{Corollary}[section]
\newtheorem{remark}{Remark}[section]

\numberwithin{equation}{section}

\begin{document}
\title[Diagonalization of the Markov Form]{Explicit Diagonalization of the Markov Form on the Temperley-Lieb Algebra}
\author{Josh Genauer}
\address{Dept. of Math., Stanford University, Standford, CA}
\email{genauer@math.stanford.edu}
\author{Neal W. Stoltzfus}
\address{ Department of Mathematics, Louisiana State University,
Baton Rouge, LA 70803-4918}
\email{stoltz@lsu.edu}

\thanks{Research partially supported by NSF Grant DMS-0353722, and a Louisiana Board of Regents Enhancement grant, LEQSF(2002-2004)-ENH-TR-17.}

\begin{abstract}
	 In fundamental paper in 1984, Vaughn Jones developed his new polynomial invariant of knots using a Markov trace on the Temperley-Lieb algebra. Subsequently, Lickorish used the associated bilinear pairing to provided an alternative proof for the existence of the 3-manifold invariants of Witten, Reshetinkin, and Turaev. A key property of this form is the non-degeneracy of this form except at the parameter values ${\pm}2{\cos}\frac{\pi}{n+1}$ \cite{Li91}.  Ko and Smolinsky  derived a recursive formula for the determinants of specific minors of Markov's form, establishing the needed non-degeneracy \cite{KoSm91}.  In this paper, we  define a triangular change of basis in which the form is diagonal and explicitly compute the diagonal entries of this matrix as  products of quotients of Chebyshev polynomials, corroborating the determinant computation of Ko and Smolinsky.  The method of proof employs a recursive method for defining the required orthogonal basis elements in the Temperley-Lieb algebra, similar in spirit to Jones'  and Wenzl's recursive formula for  a family of projectors in the Temperley-Lieb algebra.  We define a partial order on the non-crossing chord diagram basis and give an explicit formula for a recursive construction of an orthogonal basis, via a recursion over  this partial order.  Finally we relate this orthogonal basis to bases constructed using the calculus of trivalent graphs developed by Kauffman and Lins\cite{KL94}.
	
\end{abstract}
\maketitle
\begin{section}{Introduction}
In 1984, Vaughn Jones developed his new polynomial invariant of knots using a Markov trace, $t$, on the Temperley-Lieb algebra.  The Markov trace can be extended to a bilinear pairing using the algebra structure, $<x,y>=t(xy)$. In the basis given by non-crossing chord diagrams (disjoint arcs in the upper half plane with their endpoints on  $\{1,2,3,\ldots,2n\}$), the pairing can be defined geometrically by joining the first non-crossing chord diagram with the reflection of the second non-crossing chord diagram in the lower half-plane, counting the number of circles, $n$, and associating to this pair of non-crossing chord diagrams, the monomial, $q^n$, where $q$ is an indeterminate, called the \textit{quantum parameter}.  Using the theory of quantum groups and representations of lie algebras, Witten, Reshetinkin, and Turaev proved the existence of new invariants for three manifolds and links in three manifolds \cite{ReTu91,Wi89}. Using the Markov associated bilinear form,  Lickorish \cite{Li91} gave an alternative proof of the existence of these quantum invariants of a 3-manifold.
 
  A key property of this form is its non-degeneracy except when specialized at the real numbers, ${\pm}2{\cos}\frac{\pi}{n+1}$, the real parts of complex roots of unity and the roots of the Chebychev polynomials \cite{Li91}.  Ko and Smolinsky  derived a recursive formula for the determinants  of Markov's form in the non-crossing chord diagram basis and established the non-degeneracy needed by Lickorish to complete his definition\cite{KoSm91}.  
  
 In this paper, we  define a triangular change of basis to an orthogonal basis for the form and explicitly compute the diagonal entries of this matrix as  products of  quotients of Chebyshev polynomials.  
To accomplish this, we define a new partial order on the set of \textit{non-crossing} chord diagrams (which form a geometric basis for the Temperley-Lieb algebra) and use this partial order to compute an orthogonal basis, by a recursion over this poset.  

In a related setting of Topological Quantum Field Theory, Blanchet, Habegger, Masbaum and Vogel have exploited a related bilinear form and explicit orthogonal bases, see Lemma 4.1 in \cite{BHMV91} 
Pat Gilmer has also communicated the fact that other orthogonal bases for the Temperley-Lieb algebra were previously known, see the comment in \cite{GH04}, and can be constructed using the method of trivalent graphs discussed in Kauffman and Lins\cite{KL94}. In the penultimate section, we  relate our orthogonal basis to those constructed by trivalent graph techniques.

\end{section}
\begin{section}{Preliminaries}
\begin{definition}
A \textit{non-crossing chord diagram} is an isotopy class  of a set of $n$ non-intersecting properly embedded arcs   (a configuration) in the closed upper half-plane joining $2n$ points on the real line, the boundary of the upper half plane.  
\end{definition}

In this paper, we consider only non-crossing chord diagrams, as they provide a geometric basis for the Temperley-Lieb algebra. Henceforth, we drop the adjective, non-crossing. Note also that the isotopy, in the definition above, is not required to be fixed on the boundary, although it must be proper, that is, taking the boundary into itself.
 In each isotopy class, there is a unique chord diagram in \textit{standard position} defined by the condition that  the arcs are semicircles with centers on the real axis whose endpoints are located at the positive integers, $\{1,2,3,\ldots,2n\}$. We will often use this configuration for a chord diagram implicitly within our presentation.  We adopt the perspective and notation of Ko and Smolinsky \cite{KoSm91} for the set, $D_n$, of chord diagrams. 

\begin{remark}
The cardinality of the set, $D_n$ of $n$-chord diagrams  is the $n$th Catalan number, $C_n={\frac{1}{n+1}}{{2n}\choose{n}}$
\cite{Br77}.
\end{remark}

	There are many equivalent combinatorial reformulations for chord diagrams (see R. Stanley \cite{RS}, for 63 variants of the  ``Catalan disease"),  but we will mention only two relevant equivalent descriptions.  Graphically, we can uniquely associate to each chord diagram the  \textit{rooted plane tree} which is dual to the arcs (and rooted to the top). 
A second  algebraic perspective is to associate the pairing of the  end-points of a given non-intersecting chord diagram.  This is a fixed-point-free involution (\textit{perfect matching} in the combinatorial vernacular) on the set $\{1,2,3,\ldots,2n\}$  satisfisfying the condition that the orbit transpositions are non-crossing.

\end{section}\begin{section}{The Markov Symmetric Bilinear Form}


	Consider the ring $R=\mathbb{Z}[q,q^-1]$ and $M$, the free $R-$module generated by the chord diagrams.   Let $\mathbb{Q}(q)$ be the quotient field of $R$, the field rational functions (quotients of polynomials with rational coefficients), and denote by $V_n$ the free vector space spanned by $D_n$, that is, the set of formal finite sums of $n$-chord diagrams with coefficients in $\mathbb{Q}(q)$. A basis for this vector space will be denoted $e_{(a_n,\ldots,a_1)}$ and ordered by the lexicographic order.  
Following Lickorish \cite{Li91},  we define a bilinear symmetric form:
$ \langle\ \ ,\ \ \rangle_n:V_n \times V_n\to{\mathbb{Q}(q)}$
to be the bilinear extension of the  geometric pairing on the chord diagram basis given by associating the monomial, $q^c$, where $c$ is the count of the number of closed curves obtained when two chord diagrams are \textit{glued} together. 

Explicitly, we consider two $n$-chord diagram $\alpha$ and $\beta$ in standard position with endpoints at $\{1,\dots,2n\}$ on the real line $y=0$ in the plane. Reflect $\beta$ across the real axis, and then identify each node of $\alpha$ with its corresponding node in $\beta$ (the correspondence is determined by connecting endpoints with the same index). This defines a configuration we will denote $\alpha\cup\beta$. Count the number of resulting circles in $\alpha\cup\beta$, and call this number
$c:=\#(\alpha\cup\beta)$. We then define ${\langle}\alpha,\beta{\rangle}_{n}:=q^c$.  This defines the Markov form on a linearly independent basis for $V_n$ and we extend the form by bilinearity to all of ${V_n}\times{V_n}$.  The following figure illustrates the construction: 
\begin{figure}[!ht]
\includegraphics[scale=.45]{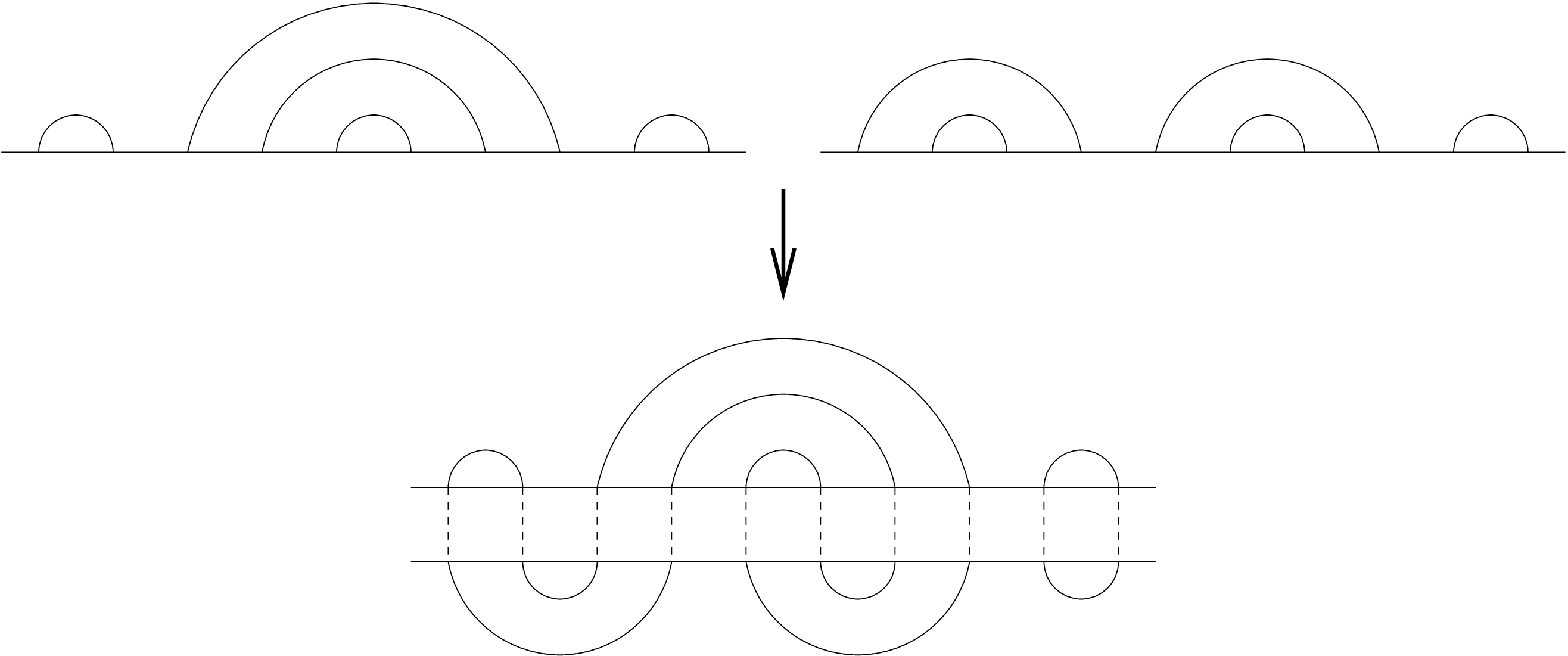}
\caption{${{\langle}(1,3,2,1,1),(2,1,2,1,1){\rangle}}=x^2$}
\end{figure}

\end{section}

\begin{section}{Chebyshev Polynomials}

Let $T_{k}$ denote the $k \times k$ tri-diagonal matrix:\\
$$\left(
\begin{matrix}
q & 1 &  &  &  & \\
1 & q & 1 &  & 0 & \\
 & 1 & q & 1 & & \\
& & \ddots & \ddots & \ddots & \\
 & 0 & & 1 & q & 1\\
 & & & & 1 & q 
\end{matrix}\right)
$$
and let the $k$th Chebyshev polynomial $\Delta_{k}$ be defined as det$(T_{k})$, a polynomial in $q$. By expanding the determinant along its top row we find the following recursive relation for for $k{\ge}3$: 

\begin{align}\label{recur}
\Delta_{k} = q{\Delta_{k-1}}-{\Delta_{k-2}}  
\end{align}

Defining $\Delta_{0}:=1$ and $\Delta_{-1}:=0$  the above identity holds for $k{\ge}1$.

\end{section}
\begin{section}{Operations on $n$-chord diagrams}

	We also need the following operations on chord diagrams from Ko and Smolinsky. For $0 < k < 2n+1$, let $l_k:D_{n}\rightarrow{D_{n+1}}$
be defined by taking a chord diagram in standard position  and adding an arc in the upper half plane whose endpoints lie on the real line between the points $(k-1,0)$ and $(k,0)$ so that the inserted arc is sufficiently small to be a new inner loop, i.e. it does not intersect any other arcs. Finally, we take the isotopy class of this new chord diagrams with $2n+2$ endpoints in $D_{n+1}$.

\begin{figure}[!ht]
\scalebox{.55}{\includegraphics{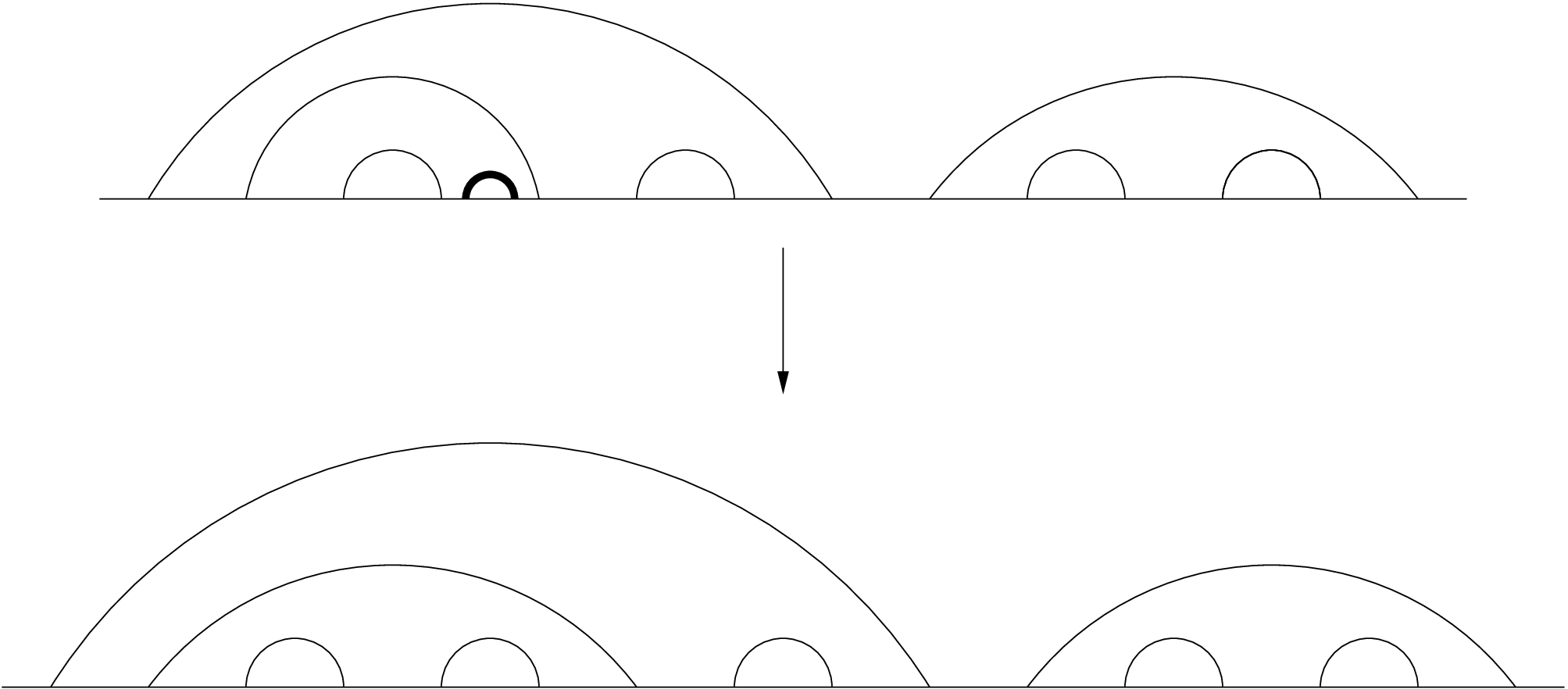}}
\caption{$l_5$}
\end{figure}

For $0<k<{2n+1}$, we also define:
${\tau{_k}}:{D_{n+1}}\rightarrow{D_{n}}$
by taking a chord diagram and adding the interval $[k,k+1]$ on the $x$-axis.  If $k$ and $k+1$ were already joined by an arc we remove the resulting component.  The result is a chord diagram on $2n-2$ vertices and we take its isotopy class in $D_{n-1}$.

\begin{figure}[!ht]
\scalebox{.65}{\includegraphics{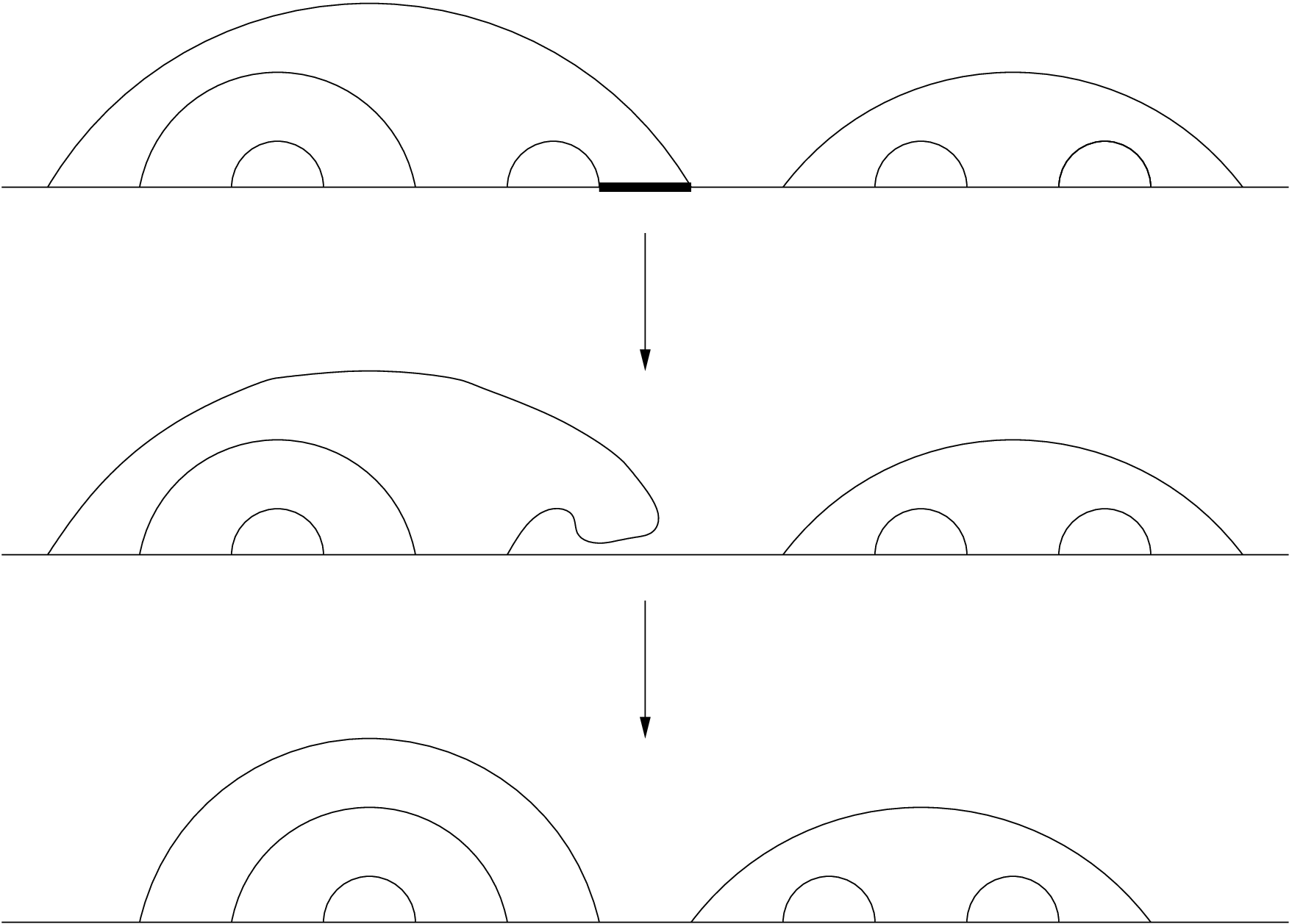}}
\caption{$\tau_7$}
\end{figure}

\begin{remark}{\label{taul}}
Note that ${\tau}_{j}{l_{k}}$ is the identity on $D_{n}$ for $j{\in}\{k-1,k,k+1\}$.
\end{remark}

\begin{remark}
The number of arcs is finite, so that  every chord diagram must have an innermost arc.
\end{remark}
Finally, we extend the operations, $l_k$ and $\tau_k$ linearly to the vector space, $V$, spanned by the chord diagrams.

	Hence any ${\alpha}\in{D_n}$  there is a ${\beta}\in{D_{n-1}}$ such that $\alpha=l_k(\beta)$ for some $k$.  Namely set $\beta={\tau_k}(\alpha)$ where $\alpha$ has an innermost arc with endpoints $k$ and $k+1$.  Moreover $k$ is uniquely determined if we insist that $k$ is as small as possible.  Recursively we obtain a unique sequence $(l_{k_n},l_{k_{n-1}},...,l_{k_{1}})$ such that $\alpha = l_{k_n}(l_{k_{n-1}}(\cdots(l_{k_{1}}(\phi))))$, where $\phi$ denotes the empty chord diagram with no arcs.  For brevity's sake we will often omit the $l$'s from the $n$-tuple and refer only to the unique sequence $(k_n,..,k_1)$. Following Ko and Smolinsky, we call this tuple the $\mathit{restricted}$ sequence associated to the chord diagram $\alpha$.

\begin{remark}
Given an arbitrary sequence $(k_n,...,k_1)$ with $0<k<{2n-1}$ one can associate the element of $D_n$ given by $l_{k_n}(l_{k_{n-1}}({\cdots}(l_{k_{1}}(\phi))))$.  This less restricted association is, however,  not one-to-one, so we shall only use restricted sequences (unless explicitly noted otherwise).
\end{remark}

\begin{figure}[!ht]
\scalebox{.55}{\includegraphics{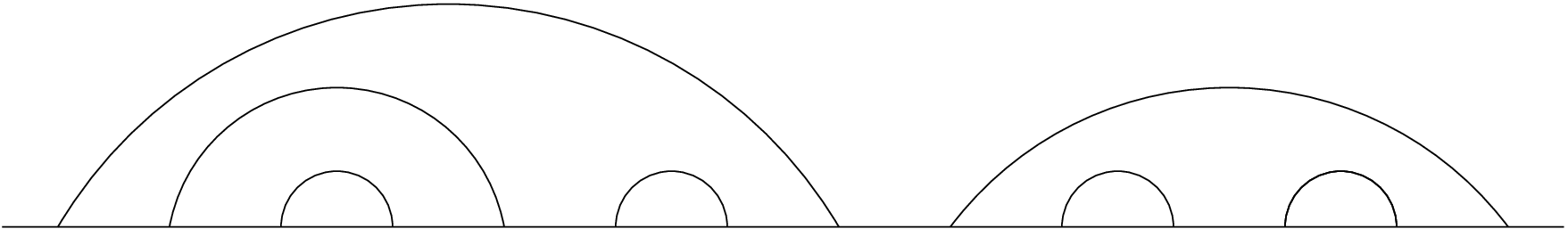}}
\caption{The chord diagram with associated restricted sequence $(3,2,2,1,2,2,1)$}
\end{figure} 

	The following simple characterization of restricted sequences will be quite useful.

\begin{proposition}
A sequence of positive integers $(k_n,...,k_1)$ is restricted if and only if $k_{1}=1$, ${k_{i+1}}\leq{k_{i}}+1$ for $1{\leq}i{\leq}n-1$. 
\end{proposition}

 First note that, given a sequence $(k_n,...,k_1)$, the leftmost inner arc has its left endpoint at $(k_n,0)$.  The only constraint when adding a new arc is to guarantee that it is the leftmost inner arc, that is, its left endpoint must  be in the range from $(1,0)$ to $({k_n}+1,0)$.  
\end{section}

pen
\begin{section}{Coordinate-wise Partial Order on Restricted Sequence}
	Using the restricted sequences, we will now define a coordinate-wise partial order on the chord diagrams.  	
\begin{definition}
	 $\alpha=(a_n,\ldots,a_1)\preccurlyeq \beta=(b_n,\ldots,b_1)=$  provided that, for all $k$ satisfying $1\leq{k}\leq{n}$, we have ${a_k}\leq{b_k}$. 
\end{definition}

The cover diagrams for these finite ranked posets are drawn below for the cases, $n=3, 4$ and $5$. 


\begin{figure}[!ht]
\scalebox{.5}{\includegraphics{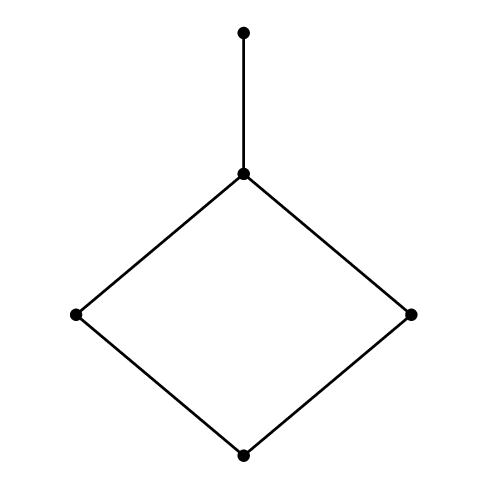}}
\caption{Cover Diagram for the Poset of Non-crossing 3-chords}
\end{figure} 

\begin{figure}[!ht]
\scalebox{.5}{\includegraphics{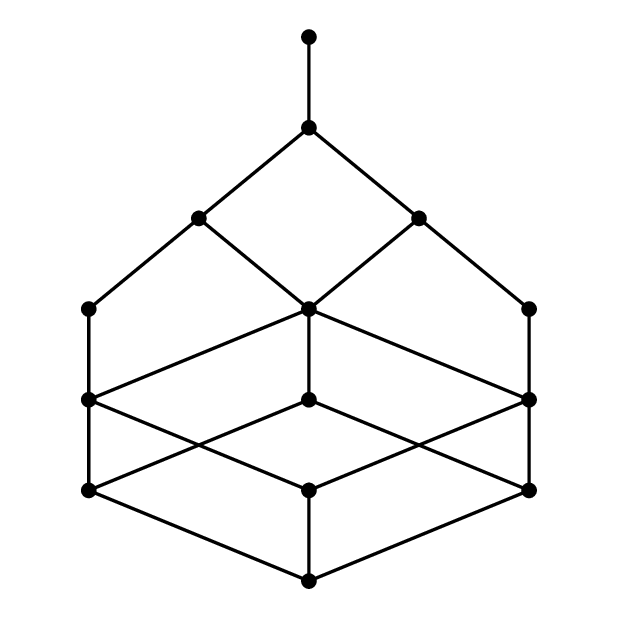}}
\caption{Cover Diagram for the Poset of Non-crossing 4-chords}
\end{figure}

\begin{figure}[!ht]
\scalebox{.5}{\includegraphics{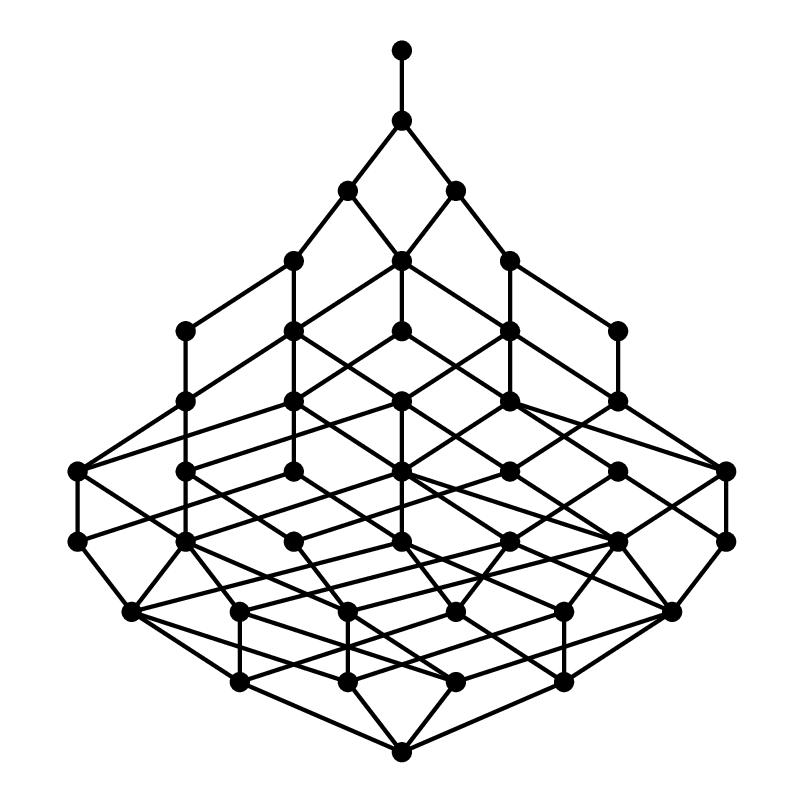}}
\caption{Cover Diagram for the Poset of Non-crossing 5-chords}
\end{figure}


Finally, we prove a proposition concerning this partial order that will be needed later.
\begin{proposition}\label{littleprop}
Let $\beta{\in}D_n$ be a chord diagram with associated restricted sequence $(b_n,...,b_1)$.  If $j{\notin}\{{b_n}-1,b_n,{b_n}+1\}$, then the restricted sequence associated to ${\tau_{j}}(\beta)$ is  $(b'_{n-1},...,b'_1)$ where:
\begin{enumerate}
\item{$b'_{n-1}={b_n}-2$ if $j<{b_n}-1$}
\item{$b'_{n-1}{\leq}{b_n}$ if $j>{b_n}+1$.}
\end{enumerate}
\end{proposition}

\begin{proof}
	For the first statement, note that the arcs whose endpoints are $j$ and $j+1$ must have the endpoints $p$ and $q$ to the right of the leftmost inner arc.  Otherwise, there would be an inner arc to the left of the leftmost inner arc.  Applying $\tau_j$, the new chord diagram is obtained from the old by deleting two arcs and joining $p$ and $q$ by an arc which does not affect the leftmost inner arc's position.  By deleting two arcs to the left of the inner arc, the inner arc has shifted two spots to the left.  Thus $b'_{n-1}={b_n}-2$.

	To verify second statement, we consider the endpoints of the arcs whose endpoints are at $j$ and $j+1$.  If there is an inner loop at $j$ (from $j$ to $j+1$), then removing that inner arc does nothing to affect the position of the leftmost inner arc.  Then $b'_{n-1}=b_n$.  Similarly if one arc or both arcs have their end points to the right of the leftmost inner arc, then $\tau_j$ does not change the position of the leftmost inner arc.  If both arcs have endpoints to the left of the leftmost inner arc, then these points must be $j'$ and $j'+1$ where $j'+1{<}b_n$ (they must be consecutive points, or else there would be an inner arc to the left of the leftmost inner arc).  But now $\tau_j$ changes the chord diagram by deleting the two arcs at $j$ and $j+1$ and adding an arc from $j'$ to $j'+1$
and leaving the rest of the diagram unchanged.  At the conclusion there is an inner arc to the left of $b_n$, so $b'_{n-1}{<}b_n$.  
\end{proof}

\end{section}
\begin{section}{Operations \& the Markov Form}

	We now state explicitly several basic facts relating the previously defined operations and the Markov form used  implicitly in \cite{KoSm91}.
\begin{proposition}{\label{fund1}}
If $\alpha,\beta{\in}D_n$, then for $1{\leq}k{\leq}2n+1$,
$\langle{{l_k}(\alpha)},{l_k}({\beta})\rangle=q{\cdot}{\langle{\alpha},{\beta}\rangle}$.
\end{proposition}
\begin{proof}
Geometrically, the inserted arcs in $l_k(\alpha)$
and $l_k(\beta)$ join to form an additional circle disjoint from
those formed by $\alpha\cup\beta$.
\begin{figure}[!ht]
\scalebox{.50}{\includegraphics{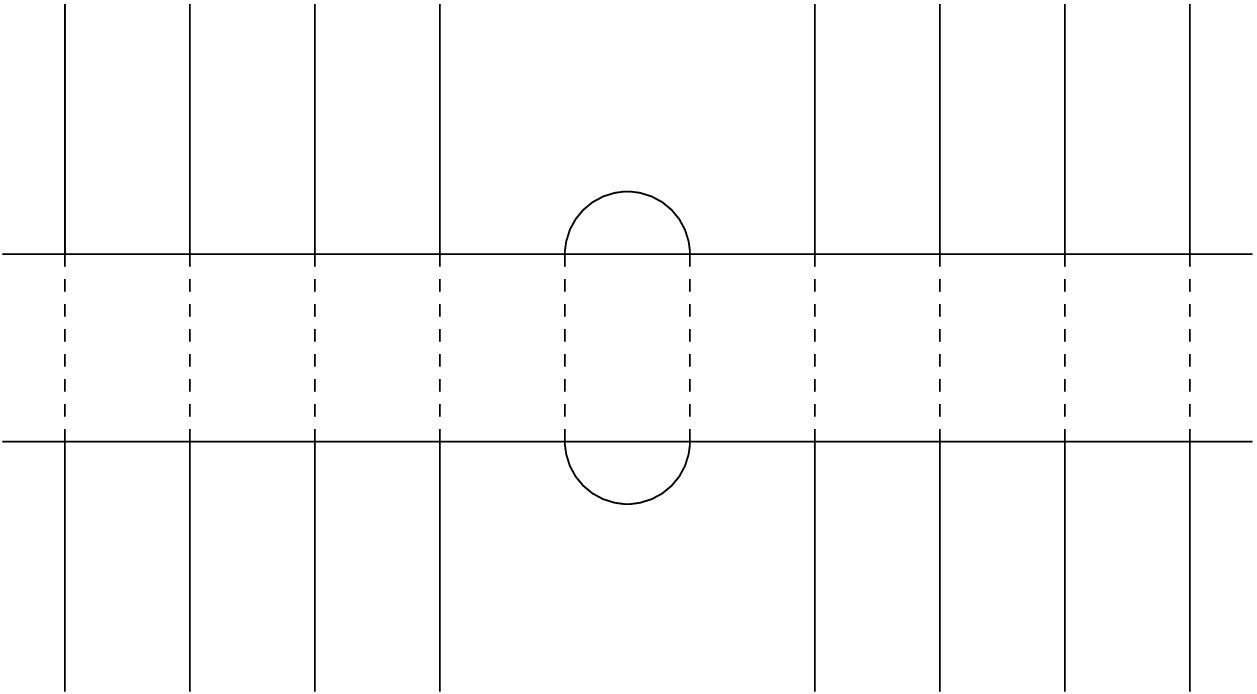}}
\caption{$\langle{{l_k}(\alpha)},{l_k}({\beta})\rangle$}
\end{figure}
\end{proof}

\begin{proposition}{\label{fund2}}
If $\alpha,\beta{\in}D_n$, then for $1{\leq}k{\leq}2n$,
$\langle{{l_{k+1}}(\alpha)},{l_k}({\beta})\rangle={\langle{\alpha},{\beta}\rangle}$
\end{proposition}
\begin{proof} Again, draw the picture.  In this case the inserted arcs form an interval which is in some component of $\alpha\cup\beta$. The number of components is unchanged.

\begin{figure}[!ht]
\scalebox{.50}{\includegraphics{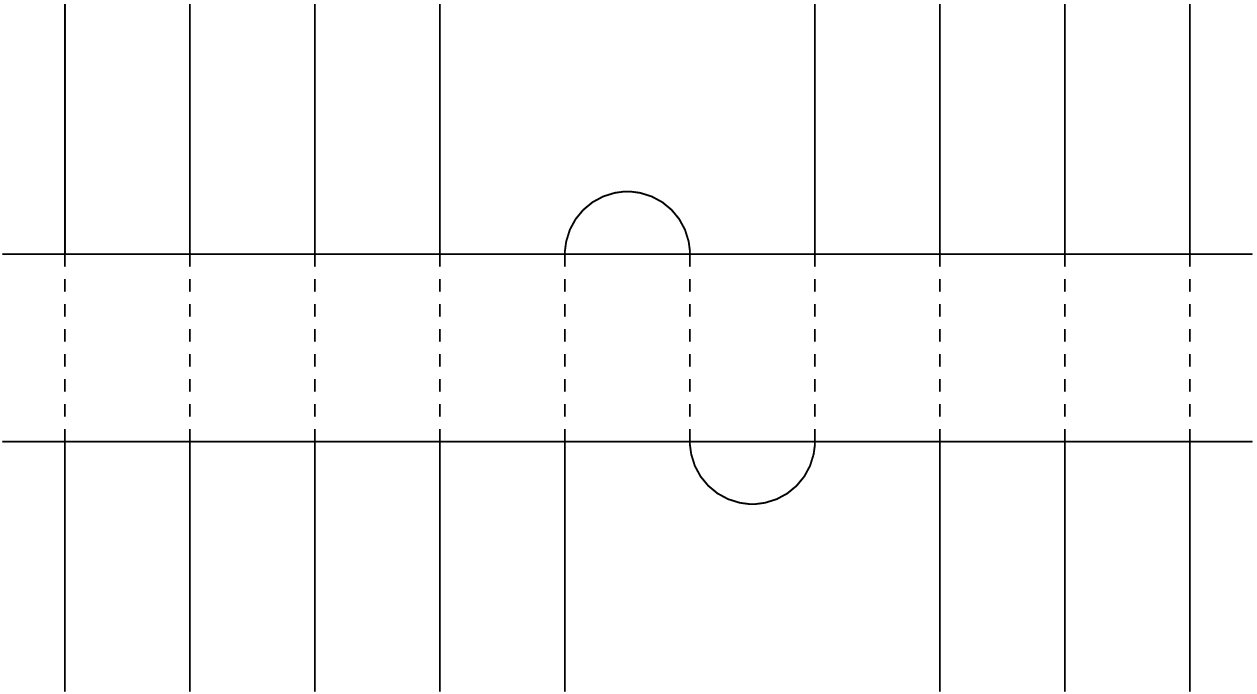}}
\caption{$\langle{{l_k}(\alpha)},{l_{k+1}}({\beta})\rangle$}
\end{figure}
\end{proof}

\bigskip

Upon linear extension to $V$,  Proposition~$\ref{fund1}$ and Proposition~$\ref{fund2}$ imply:
\begin{proposition}\label{extend}
For $f,g{\in}V_n$ and $1{\leq}k{\leq}2n+1$, 
\begin{eqnarray*}
q{\cdot}{{\langle}{f},{g}{\rangle}}=q{\cdot}{\langle}{{l_k}(f)},{{l_{k+1}}(g)}{\rangle}={\langle}{{l_k}(f)},{{l_k}(g)}{\rangle}
\end{eqnarray*}
\end{proposition}
\begin{proof}
The statement holds since $l_k$ is linear and the conditions are linear in $j$ and $g$.
\end{proof}

\begin{proposition}
If $\alpha{\in}D_{n+1}$ and $\beta{\in}D_n$ then for $1{\leq}k{\leq}2n$, ${q^c}{\langle}{\tau}_{k}(\alpha),{\beta}{\rangle}_{n}={\langle}({\alpha}),{l_{k}}({\beta}){\rangle}_{n+1}$ where $c=1$ if $\alpha$ has an inner loop from $k$ to $k+1$ and $c=0$ otherwise.
\end{proposition}
\begin{proof}
Consider the join ${\alpha}{\cup}{l_{k}}({\beta})y$. $l_k(\beta)$ has an inner loop from $k$ to $k+1$ so the lower part of the gluing at the interval $[k,k+1]$ resembles the letter \textbf{U}, with its bottom half below the x-axis.  Push the \textbf{U} onto the x-axis obtaining the figure $\sqcup$ with its base resting along the x-axis.  The resulting chord diagram is none other than ${\tau}_{k}(\alpha){\cup}{\beta}$ unless $\beta$ originally had an arc from $k$ to $k+1$.  In that case there is an extra loop in the right expression.
\end{proof}

\end{section}\begin{section}{Construction of an Orthogonal Basis for $V_n$}
	We now begin the construction our desired orthogonal basis elements using a two-level recursion.  The outer recursion is on $n$ itself, the size of the chord diagrams, the inner recursion is on the coordinate-wise partial ordering of restricted sequences associated to the chord diagrams.  We now define our recursive expression for calculating the orthogonal basis elements (for $2{\leq}a_n{\leq}n$):

\begin{align}
e'_{({a_n},...,{a_1})}=l{_{a_n}}(e'_{({a_{n-1}},...,{a_1})})-\frac{\Delta_{{a_n}-2}}{\Delta_{{a_n}-1}}{\cdot}{e'_{({{a_n}-1},{a_{n-1}},...,{a_1})}}
\end{align}
Note that this definition requires a recursion on $n$ as well as the coordinate-wise
partial order, $\preccurlyeq$, on $D_n$. 

We initialize our recursion with the definition:

\begin{eqnarray*}
e'_{(1)}:=e_{(1)}
\end{eqnarray*}

With the convention $\Delta_{-1}=0$, the recursive formula makes sense when
$a_n=1$ (although $(0,a_{n-1},\ldots,a_1)$ is an invalid restricted sequence) and
simplifies to:
\begin{eqnarray*}
e'_{(1,{a_{n-1}},\ldots,{a_1})}=l_{1}(e'_{({a_{n-1}},\ldots,{a_1})})
\end{eqnarray*}
which is a recursion solely on $n$.

\end{section}

\begin{remark} 
We can use the map $l_1$ to embed $V_{n-1}$ into $V_n$
(this sends the restricted sequence $(a_{n-1},a_{n-2},\ldots,a_1)$ to $(1,a_{n-1},a_{n-2},\ldots,a_1)$. The Markov form is scaled by $q$ arising from the one additional component (by Prop~$\ref{fund1}$).
\end{remark}

\begin{section}{Quad moves}

In this section, we investigate the relationship between the operation $l_k$ and the partial ordering $\preccurlyeq$.

\begin{lemma}{\label{prec}}
Let ${\alpha}$ and ${\beta}{\in}D_n$.  If $1{\leq}{k}{\leq}2n+1$ then:
\begin{align}
\alpha{\preccurlyeq}{\beta}{\Leftrightarrow}l_k(\alpha){\preccurlyeq}l_k(\beta)
\end{align}
\end{lemma}

The proof, given at the end of the section, will depend on an operation on chord diagrams called the \textit{quad moves} --- equivalently \textit{0-surgery} on a chord diagram.

\begin{definition} For each pair of standard arcs, $\overline{ik}$ and $\overline{jl}$ in a chord diagram $\alpha$ with endpoints satisfying $i<k<j<l$ and such that no arc with endpoint in the interval $[k,j]$ separates the given arc, we define a quad move, $Quad_{i,j}(\alpha)$, to be the chord diagram obtained from $\alpha$ by removing the given pair and replacing them by the arcs $\overline{jk}$ and $\overline{il}$.
\end{definition}

In the dual \textit{rooted plane tree}, the condition on the arcs requires the corresponding edges to have a common parent node.

	There is an important relationship between quad moves and the partial ordering.  First, we demonstrate that the result of a  quad move on a diagram is smaller in the partial order, $\preccurlyeq$, than the original diagram.  Furthermore, we show for $\alpha{\preccurlyeq}\beta$ in $D_n$ there is a sequence of quad moves which transforms $\beta$ into $\alpha$.

\begin{figure}[!ht]
\scalebox{.65}{\includegraphics{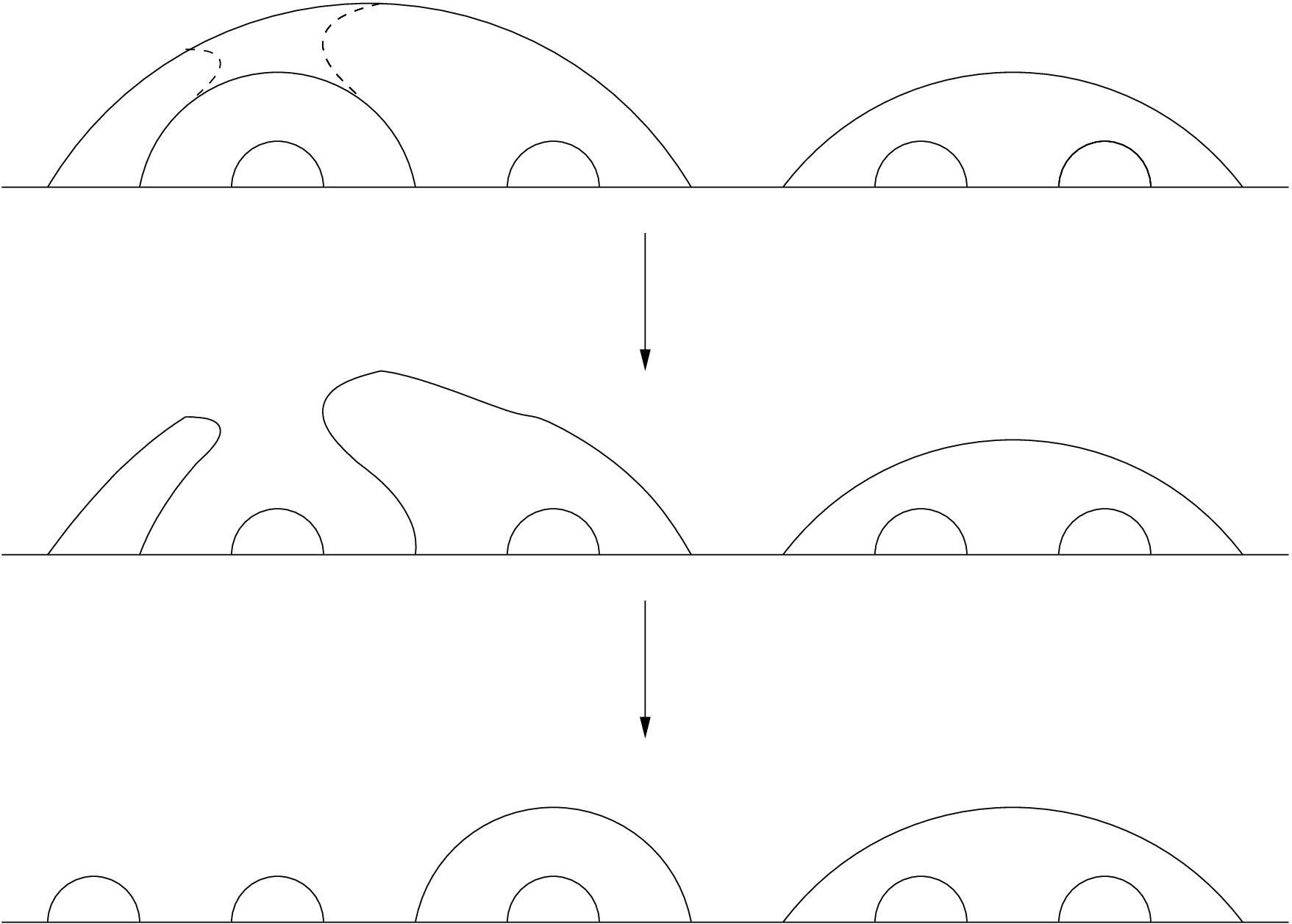}}
\end{figure}
\begin{definition}
The \textit{depth} of an arc $a$ in $\alpha{\in}D_n$ is the number of arcs that lie over it. An arc $a$ lies over another arc $b$ if the left endpoint of $a$ lies to the left of left endpoint $b$ and the right endpoint of $b$ lies to the right of the right endpoint of $b$.
\end{definition}

\begin{remark}
The depth of the $i$th arc in ${\alpha}{\in}D_n$ is $a_i$ where $(a_n,...,a_1)$ is the restricted sequence of $\alpha$.  Additionally, from the \textit{rooted plane tree} perspective, depth is the distance to the root.
\end{remark}

\begin{proposition}
If $(a_n,...,a_1)=\alpha{\in}D_n$ and $Q_{k,k'}$ a quad move of $\alpha$ involving arcs $a_k$ and $a_k'$ with $k{<}k'$, then $\beta=Q_{k,k'}(\alpha){\preccurlyeq}\alpha$.
\end{proposition}

\begin{proof}
We show that the depth of the $i$th arc of $\beta=(b_n,...,b_1)$ is less than or equal to depth of the $i$th arc of $\alpha$.  The depth of 
the $j$th arc for $j{<}{k'}$ remains unchanged since they lie entirely to the right of the arcs involved, so $b_j=a_j$.  Let $\{a_{k'+1},...,a_{k''}\}$ be the
arcs that lie under $a_{k'}$.  Then clearly for ${k'+1}{<}j{<}{k''}$ we must have $a_{k'}{<}a_j$ and also that $b_{j-1}=a_{j}-2$.  But this implies 
that $b_{j-1}=a_{j}-2{<}{a_{j}-1}{\leq}{a_{j-1}}$. Thus, for $1{\leq}j{\leq}{k''-1}$ we have $b_j{\leq}a_j$.  Now the depth of the $k''$th arc of $\beta$
is one less than the depth of the $k'$th arc of $\alpha$, so $b_{k''}={a_{k'}}-1$.  But then $b_{k''}={a_{k'}}-1{<}{a_{k'}}{<}{a_{k''}}$.  The depths of 
the other arcs remain unchanged under the quad move since they lie entirely to the left of it.
\end{proof}

\begin{proposition}
Let $\alpha$ and $\beta$ be in $D_n$.  If $\alpha{\preccurlyeq}\beta$, then there is a 
sequence of quad moves which transforms $\beta$ into $\alpha$.
\end{proposition}
\begin{proof}
Note that removing an inner arc from a chord diagram does not affect the depth of the other arcs.  In terms of restricted sequences, the equivalent way to express the deletion of an inner arc, which say is the $i$th arc in the diagram, would be to delete the $i$th component of the restricted sequence where the $i$th arc is an inner arc and shift all components to the
left of the $i$th component to the right one space, thus obtaining an $n-1$ restricted sequence.
By the condition $\alpha{\preccurlyeq}\beta$, $\alpha$'s rightmost inner arc lies to the right of $\beta$'s.  Performing quad moves on $\beta$ involving its rightmost inner arc and the arc immediately to its right, we may \textit{pass the inner arc under} however many arcs on its right, until it aligns with the inner arc of $\beta$.  This quad move leaves the depth of all arcs
unaltered except for the depth of the right most inner arc which is lowered by one.  The quad move leaves the ordering of the arcs unchanged except for interchanging the rightmost inner arc with the arc immediately to its left.  Thus the quad move of $\beta$ is still $\preccurlyeq{\alpha}$.  So now we can
assume that the rightmost inner arc of $\beta$ is aligned with the rightmost inner arc of $\alpha$.  Let ${\alpha}',{\beta}'{\in}D_{n-1}$ be the chord
diagrams obtained by removing the rightmost inner arcs from $\alpha$ and $\beta$ respectively.  Because the rightmost inner arcs are in the same position, we must have ${{\alpha}'}{\preccurlyeq}{{\beta}'}$.  Then by induction on the size of chord diagrams there is a sequence of quad moves which take ${{\beta}'}$ to ${{\alpha}'}$.   (The induction hypothesis is trivially satisfied for one chord diagrams.) Combining the two sequences of quad moves completes the proof.
\end{proof}

\begin{proof}
To prove the initial lemma, we need only verify that quad moves commute with inserting inner arcs and quad moves also commute with deleting inner arcs, as long as
the inner arc is not the result of a quad move.
\end{proof}

We now give two corollaries of the relationship between the partial order and the $l_k$ operation.

\begin{corollary}
The new basis $e'_{(a_n,...,a_1)}$ for the space, $V_n$ is a linear combination of the original basis elements such that the coefficient of $e_{(b_n,...,b_1)}$ is nonzero in the linear combination iff $(b_n,...,b_1){\preccurlyeq}(a_n,...,a_1)$.
\end{corollary}
\begin{proof}
$e'_{({a_n},...,{a_1})}=l{_{a_n}}(e'_{({a_{n-1}},...,{a_1})})-\frac{\Delta_{{a_n}-2}}{\Delta_{{a_n}-1}}{\cdot}{e'_{({{a_n}-1},{a_{n-1}},...,{a_1})}}$.
It suffices to show that $e'_{({{a_n}-1},{a_{n-1}},...,{a_1})}$ is a linear combination of of basis elements in the \textit{downset} of $\alpha=({a_n},\ldots,{a_1})$, the set $D(\alpha):=\{\beta|\beta\preccurlyeq \alpha\}$. Inductively $e'_{({{a_n}-1},{a_{n-1}},\ldots,{a_1})}$ is a linear combination of of basis elements $D({{a_n}-1},\ldots,{a_1})$.  But by the transitivity of the partial ordering, this implies these basis elements are in $D({a_n},\ldots,{a_1})$.  For the other term, we may inductively assume that $e'_{({a_{n-1}},\ldots,{a_1})}$ is a linear combination of elements all $D({a_{n-1}},\ldots,{a_1})$.  Noticing that $l{_{a_n}}(e_{({a_{n-1}},\ldots,{a_1})})=e_{({a_n},\ldots,{a_1})}$ we conclude from the preceding lemma, Lemma~$\ref{prec}$ that for any basis element $e$ such that $e{\preccurlyeq}e_{({a_{n-1}},\ldots,{a_1})}$ we have $l{_{a_n}}(e){\preccurlyeq}{e_{({a_n},\ldots,{a_1})}}$. The result follows.
\end{proof}

Finally we determine the coefficient of the largest term (for $\preccurlyeq $) for the new basis.

\begin{corollary} The change of basis is triangular, that is, 
in the expression of $e'_{(a_n,\ldots,a_1)}$ as a linear combination of elements of the original basis , the polynomial coefficient of $e_{({a_n,\ldots,a_1})}$ is one.
\end{corollary}
\begin{proof}
Note that only $l{_{a_n}}(e'_{({a_{n-1}},\ldots,{a_1})}$ has a nonzero $e_{({a_{n}},\ldots,{a_1})}$ component.  Inductively we assume that $e'_{(a_{n-1},\ldots,a_1)}$ written as a linear combination of the original basis elements has 1 for the coefficient of $e_{({a_{n-1},\ldots,a_1})}$. By the conditions for equality in Lemma~$\ref{prec}$, the result follows.
\end{proof}

\end{section}

\begin{section}{Orthogonality}

\begin{theorem}{\label{Fact Two}}
The following two results hold:\\

\begin{itemize}
\item $(b_n,...,b_1){<}(a_n,...,a_1)\Rightarrow{\langle{e_{(b_n,...,b_1)}},{e'_{(a_n,...,a_1)}}\rangle}=0$
\item $(b_n,...,b_1)=(a_n,...,a_1)\Rightarrow{\langle{e_{(a_n,...,a_1)}},{e'_{(a_n,...,a_1 )}}\rangle}=\frac{\Delta_{a_n}}{\Delta_{{a_n}-1}}{\cdot}{\cdot}{\cdot}\frac{\Delta_{a_1}}{\Delta_{{a_1}-1}}$
\end{itemize}
\end{theorem}
\begin{proof}
Our proof is by a double induction, an outer induction on the number of arcs, $n$, and an inner induction on the ordering of restricted sequences of length $n$. 

We first note that the result holds for $n$=1, $\langle
e_1,e'_1\rangle=q=\frac{\Delta_{1}}{\Delta_{{1}-1}}$.

 Our base case for the induction on the ordering is:
Let $(b_n,...,b_1){\leq}(a_n,...,a_1)$ with ${a_n}=1$.  Then we must have ${b_n}=1$ and\\
$(b_{n-1},...,b_1){\leq}(a_{n-1},...,a_1)$ with equality iff we have equality in the original restricted $n$-sequences.  Applying the recursion formula, we obtain:

\begin{align*}
{\langle{e_{({b_n},...,{b_1})}},{e'_{({a_n},...,{a_1})}}\rangle}&={\langle{e_{({1,b_{n-1},...,{b_1})}},{e'_{({1,a_{n-1},...,{a_1})}}}}\rangle}\\
&={\langle{e_{(1,{b_{n-1}},...,{b_1})}},l_{1}(e'_{({a_{n-1}},...,{a_1})})\rangle}\\
&-{\frac{\Delta_{-1}}{\Delta_{0}}}{\langle{e_{(1,{b_{n-1}},...,{b_1})}},{e'_{({{a_n}-1},{a_{n-1}},...,{a_1})}}\rangle}\\
\intertext{the right term is zero since $\Delta_{-1}:=0$}
&=q{\cdot}{\langle{\tau_{1}}(e_{({1,b_{n-1},...,{b_1}})}),{e'_{(a_{n-1},...,{a_1})}}\rangle}\\
&=\frac{\Delta_{1}}{\Delta_{0}}{\cdot}{\langle{e_{({b_{n-1}},...,{b_1})}},{e'_{(a_{n-1},...,{a_1})}}\rangle}.
\end{align*}
By the adjoint relation between the operations and the Markov form, we have completed the proof in our base case.

The induction step is broken into four cases.  While the calculations are different for each case, they are of similar flavor.  For the sake of brevity, we give a complete calculation for a relatively complicated case, \textbf{Case 2(ii)} and sketch the others. 

\textbf{Case 1:} ${b_n}<{{a_n}-1}$:
\begin{align*}
{\langle{e_{(b_n,\ldots,b_1)}},l_{a_n}(e'_{(a_{n-1},\ldots,a_1)})\rangle}&=(q^{c}){\langle{{\tau_{a_n}}(e_{(b_n,\ldots,b_1)})},e'_{(a_{n-1},\dots,a_1)}\rangle}\\
\intertext{which by the second part of
Proposition~$\ref{littleprop}$}
&=(q^{c}){\langle{(e_{({b'_{n-1}},\ldots,{b'_1})})},e'_{({a_{n-1}},\ldots,{a_1})}\rangle}\\
\intertext{for some restricted sequence
$({b'_{n-1}},\ldots,{b'_1})$ with $b'_{n-1}{\leq}b_n$.  And since
$b'_{n-1}{\leq}b_n{<}{{a_n}-1}{\leq}{a_{n-1}}$ by induction on
$n$,}
&=(q^{c}){\langle{e_{({b'_{n-1}},\dots,{b'_1})}},e'_{({a_{n-1}},\ldots,{a_1})}\rangle}=0.
\end{align*}

The monomial coefficient $c$ is non-trivial in the case of an inner arc in $(b_n,\ldots,b_1)$ at $a_n$.

For the rest of the cases, we use the following lemma:
\begin{lemma}{\label{sub}}
Let $(b_n,\ldots,b_1),(a_{n-1},\ldots,a_1)$ be restricted sequences.  If $j{\leq}{b_n}-2<{a_{n-1}}$ then ${\langle{e_{(b_n,\ldots,b_1)}},{e'_{(j,\ldots,a_1)}}\rangle}=0$.
\end{lemma}
\begin{proof}
Induction on $j$.
For $j=1$,
\begin{align*}
{\langle{e_{(b_n,\ldots,b_1)}},{e'_{(j,\ldots,a_1)}}\rangle}&={\langle{e_{(b_n,\ldots,b_1)}},{l_1}(e'_{(a_{n-1},\ldots,a_1)})\rangle}\\
&={\langle{{\tau_1}(e_{(b_n,\ldots,b_1)})},e'_{(a_{n-1},\ldots,a_1)}\rangle}\\
\intertext{By the first part of Proposition~$\ref{littleprop}$, ${\tau_1}(b_n,\ldots,b_1)=(b'_{n-1},\ldots,b'_1)$ where $b'_{n-1}={{b_n}-2}<{a_{n-1}}$.  By the induction hypothesis on $n$,}
&={\langle{\tau_{1}}(e_{(b_n,\ldots,b_1)}),e'_{(a_{n-1},\ldots,a_1)}\rangle}\\
&=0.
\end{align*}
For $1<j{\leq}{b_n}-2$,
\begin{align*}
{\langle{e_{({b_n},\ldots,{b_1})}},{e'_{(j,\ldots,{a_1})}}\rangle}&={\langle{{\tau_j}(e_{({b_n},\ldots,{b_1}))}},e'_{({a_{n-1}},\ldots,{a_1})}\rangle}-{\frac{\Delta_{j-2}}{\Delta_{j-1}}}{\cdot}{\langle{e_{({b_n},\ldots,{b_1})}},{e'_{({j-1},{a_{n-1}},\ldots,{a_1})}}\rangle}\\
\intertext{By induction on $j$, $\langle{e_{({b_n},\ldots,{b_1})}},{e'_{({j-1},{a_{n-1}},\ldots,{a_1})}\rangle}=0$.  Once again, finding that ${\tau_j}(e_{({b_n},\ldots,{b_1})})=(b'_{n-1},\ldots,b'_1)$ with $b'_{n-1}{\leq}{b_n}-2$ we conclude by the induction hypothesis on $n$ that the first term is also zero.}
\end{align*}
\end{proof}

\textbf{Case 2:} $b_n={{a_n}-1}$
 
\begin{align*}
{\langle{e_{(b_n,\ldots,b_1)}},{e'_{(a_n,\ldots,a_1)}}\rangle}&={\langle{e_{(b_n,\ldots,b_1)}},l_{a_n}(e'_{(a_{n-1},\ldots,a_1)})\rangle}\\
&-{\frac{\Delta_{{a_n}-2}}{\Delta_{{a_n}-1}}}{\langle{e_{(b_n,\ldots,b_1)}},{e'_{({a_n}-1,a_{n-1},\ldots,a_1)}}\rangle}\\
&={\langle{e_{({b_{n-1}},\ldots,{b_1})}},e'_{({a_{n-1}},\ldots,{a_1})}\rangle}\\
&-{\frac{\Delta_{{a_n}-2}}{\Delta_{{a_n}-1}}}{\langle{e_{({b_n},\ldots,{b_1})}},{e'_{({{a_n}-1},{a_{n-1}},\ldots,{a_1})}}\rangle}
\end{align*}

\textbf{Case 2($i$):} $(b_n,\ldots,b_1)=({a_n}-1,a_{n-1},\ldots,a_1)$
By the induction hypothesis we have:
\begin{align*}
&=(\frac{\Delta_{a_{n-1}}}{\Delta_{a_{n-1}-1}}{\cdot}{\frac{\Delta_{a_{n-2}}}{\Delta_{a_{n-2}-1}}}{\cdots}\frac{\Delta_{a_{1}}}{\Delta_{{a_1}-1}})
-\frac{\Delta_{{a_n}-2}}{\Delta_{{a_n}-1}}{\cdot}(\frac{\Delta_{{a_n}-1}}{\Delta_{{a_n}-2}}{\cdot}\frac{\Delta_{a_{n-1}}}{\Delta_{{a_{n-1}-1}}}{\cdot}\frac{\Delta_{a_{n-2}}}{\Delta_{{a_{n-2}-1}}}{\cdots}\frac{\Delta_{a_1}}{\Delta_{{a_1}-1}})\\
&=0.
\end{align*}

\textbf{Case 2($ii$):} ${b_n}={{a_n}-1}$ and $(b_{n-1},\ldots,b_1){\neq}(a_{n-1},\ldots,a_1)$

Expanding $e'_{(a_n,\ldots,a_1)}$ three times with the recursive definition, we obtain:
\begin{align*}
{\langle{e_{({b_n},\ldots,{b_1})}},{e'_{({a_n},\ldots,{a_1})}}\rangle}&={\langle{e_{({b_n},\ldots,{b_1})}},l_{a_n}({e'_{({a_{n-1}},\ldots,{a_1})}})\rangle}\\
&-{\frac{\Delta_{{a_n}-2}}{\Delta_{{a_n}-1}}}{\langle{e_{({b_n},\ldots,{b_1})}},l_{{a_n}-1}({e'_{({a_{n-1}},\ldots,{a_1})}})\rangle}\\
&+{\frac{\Delta_{{a_n}-3}}{\Delta_{{a_n}-1}}}{\langle{e_{({b_n},\ldots,{b_1})}},l_{{a_n}-2}({e'_{({a_{n-1}},\ldots,{a_1})}})\rangle}\\
&-{\frac{\Delta_{{a_n}-4}}{\Delta_{{a_n}-1}}}{\langle{e_{({b_n},\ldots,{b_1})}},{e'_{({a_n}-3,{a_{n-1}},\ldots,{a_1})}}\rangle}\\
\intertext{and by the adjoint property of the operations in  Proposition~$\ref{extend}$ we get}\\
&={\langle{{\tau_{a_n}}(e_{({b_n},\ldots,{b_1})})},{e'_{({a_{n-1}},\ldots,{a_1})}}\rangle}\\
&-q{\cdot}{\frac{\Delta_{{a_n}-2}}{\Delta_{{a_n}-1}}}{\langle{{\tau_{{a_n}-1}}(e_{({b_n},\ldots,{b_1})})},{e'_{({a_{n-1}},\ldots,{a_1})}}\rangle}\\
&+{\frac{\Delta_{{a_n}-3}}{\Delta_{{a_n}-1}}}{\langle{{\tau_{{a_n}-2}}(e_{({b_n},\ldots,{b_1})})},{e'_{({a_{n-1}},\ldots,{a_1})}}\rangle}\\
&-{\frac{\Delta_{{a_n}-4}}{\Delta_{{a_n}-1}}}{\langle{e_{({b_n},\ldots,{b_1})}},{e'_{({a_n}-3,{a_{n-1}},\ldots,{a_1})}}\rangle}\\
\intertext{since ${a_n}-2,{a_n}-1,a_n$ equal ${b_n}-1,{b_n},{b_n}+1$ respectively, we use Remark~$\ref{taul}$ to get}\\
&=\big{[}(1-{q}\frac{\Delta_{{a_n}-2}}{\Delta_{{a_n}-1}}+\frac{\Delta_{{a_n}-3}}{\Delta_{{a_n}-1}}){\cdot}{\langle{e_{({b_{n-1}},\ldots,{b_1})}},{e'_{({a_{n-1}},\ldots,{a_1})}}\rangle}\big{]}\\
&-{\big{[}}{\frac{\Delta_{{a_n}-4}}{\Delta_{{a_n}-1}}}{\cdot}{\langle{e_{({b_n},\ldots,{b_1})}},e'_{({{a_n}-3},{a_{n-1}},\ldots,{a_1})})\rangle}{\big{]}}.
\intertext{By the recursion property for the Chebyshev polynomials $\ref{recur}$, $(1-{q}{\frac{\Delta_{{a_n}-2}}{\Delta_{{a_n}-1}}}+\frac{\Delta_{{a_n}-3}}{\Delta_{{a_n}-1}})=1-1=0$}\\
\intertext{For $j={{a_n}-3}$ we have $j={{a_n}-3}={{b_n}+1}-3={{b_n}-2}<{b_n}={a_n}-1{\leq}{a_{n-1}}$ hence the required conditions hold for applying Lemma~$\ref{sub}$ and we conclude that ${\langle{e_{({b_n},\ldots,{b_1})}},e'_{({{a_{n}}-3}=j,{a_{n-1}},\ldots,{a_1})}\rangle}=0$.}
\end{align*}

Some extra care should be taken when $a_n$ is small.  Because of the hypothesis, ${b_n}={{a_n}-1}$, we may assume $2{\leq}a_n$.  When $4{\leq}a_n$ expanding by recursion is well defined.  For $3=a_n$ the expansion also holds, but there is no second term.  The case $2=a_n$ deserves a special treatment.  Given $(1,b_{n-1},...,b_1)$:

\begin{align*}
{\langle{e_{(1,\ldots,b_1)}},{e'_{(2,\ldots,a_1)}}\rangle}&={\langle{l_1(e_{(b_{n-1},\ldots,b_1)})},l_{2}(e'_{(a_{n-1},\ldots,a_1)})\rangle}\\
&-{\frac{\Delta_{0}}{\Delta_{1}}}{\langle{e_{(1,b_{n-1},\ldots,b_1)}},{e'_{({1},{a_{n-1}},\ldots,{a_1})}}\rangle}\\
&={\langle{l_1(e_{(b_{n-1},\ldots,b_1))}},l_{2}(e'_{(a_{n-1},\ldots,a_1)})\rangle}\\
&-{\frac{\Delta_{0}}{\Delta_{1}}}{\langle{l_1(e_{(b_{n-1},\ldots,b_1)})},{l_1(e'_{({a_{n-1}},\ldots,{a_1})})}\rangle}\\
\intertext{And by Proposition~$\ref{extend}$ we get}\\
&={\langle{e_{(b_{n-1},\ldots,b_1)}},e'_{(a_{n-1},\ldots,a_1)}\rangle}\\
-{\frac{1}{q}}{\cdot}{\langle{e_{(b_{n-1},\ldots,b_1)}},{e'_{({a_{n-1}},\ldots,{a_1})}}\rangle}\\
&=0.
\end{align*}\\

\textbf{Case 3:} $b_n=a_n$ and $(b_{n-1},\ldots,b_1){<}(a_{n-1},\ldots,a_1)$\\

When $1=a_n$ the following holds by Proposition~$\ref{extend}$ and by induction on $n$: 
\begin{align*}
{\langle{e_{(b_n,\ldots,b_1)}},{e'_{(a_n,\ldots,a_1)}}\rangle}&={\langle{l_{1}(e_{(b_{n-1},\ldots,b_1)})},l_{1}(e'_{(a_{n-1},...,a_1)})\rangle}\\
&=q{\cdot}{\langle{e_{(b_{n-1},\ldots,b_1)}},{e'_{(a_{n-1},\ldots,a_1)}}\rangle}\\
&=0
\end{align*}

For $2=a_n$ we expand twice and apply induction on $n$.

\begin{align*}
{\langle{e_{(b_n,\ldots,b_1)}},{e'_{(a_n,\ldots,a_1)}}\rangle}&={\langle{l_{2}(e_{(b_{n-1},\ldots,b_1)})},l_{2}(e'_{(a_{n-1},\ldots,a_1)})\rangle}\\
&-{\frac{\Delta_{0}}{\Delta_{1}}}{\langle{l_{2}(e_{(b_{n-1},\ldots,b_1)})},{l_{1}(e'_{({{a_n}-1},{a_{n-1}},\ldots,{a_1})})}\rangle}\\
&=q{\cdot}{\langle{e_{(b_{n-1},\ldots,b_1)}},{e'_{(a_{n-1},\ldots,a_1)}}\rangle}-{\frac{1}{q}}{\langle{e_{(b_{n-1},\ldots,b_1)}},{e'_{(a_{n-1},\ldots,a_1)}}\rangle}\\
&=0
\end{align*}

Lastly for $2{\leq}a_n$ we expand twice obtaining an expression similar to the one above:
\begin{align*}
&={\langle{l_{b_n=a_n}(e_{(b_{n-1},\ldots,b_1)})},l_{a_n}(e'_{(a_{n-1},\ldots,a_1)})\rangle}\\
&-{\frac{\Delta_{{a_n}-2}}{\Delta_{{a_n}-1}}}{\cdot}{\langle{l_{b_n=a_n}(e_{(b_{n-1},\ldots,b_1)})},{l_{{a_n}-1}(e'_{({{a_n}-1},{a_{n-1}},\ldots,{a_1})})}\rangle}\\
&+{\frac{\Delta_{{a_n}-3}}{\Delta_{{a_n}-1}}}{\langle{{e_{(b_{n},\ldots,b_1)}},e'_{({{a_n}-2},{a_{n-1}},\ldots,{a_1})}}\rangle}\\
&=(q-{\frac{\Delta_{{a_n}-2}}{\Delta_{{a_n}-1}}}){\cdot}{\langle{e_{(b_{n-1},\ldots,b_1)}},e'_{(a_{n-1},\ldots,a_1)}\rangle}\\
&+{\frac{\Delta_{{a_n}-3}}{\Delta_{{a_n}-1}}}{\cdot}{\langle{{e_{(b_{n},\ldots,b_1)}},e'_{({{a_n}-2},{a_{n-1}},...,{a_1})}}\rangle}\\
\intertext{Note, by induction on $n$, that ${\langle{e_{(b_{n-1},\ldots,b_1)}},e'_{(a_{n-1},\ldots,a_1)}\rangle}=0$}\\
&={\frac{\Delta_{{a_n}-3}}{\Delta_{{a_n}-1}}}{\cdot}{\langle{{e_{(b_{n},\ldots,b_1)}},e'_{({{a_n}-2},{a_{n-1}},\ldots,{a_1})}}\rangle}\\ 
\intertext{For $j={{a_n}-2}$ we have $j={{a_n}-2}={{b_n}-2}<{{b_n}-1}={{a_n}-1}{\leq}{a_{n-1}}$ hence we can conclude from Lemma~$\ref{sub}$ that}
&={\langle{{e_{(b_{n},\ldots,b_1)}},e'_{({{a_n}-2},{a_{n-1}},\ldots,{a_1})}}\rangle}\\ 
&=0
\end{align*}

\textbf{Case 4:} $(b_n,\ldots,b_1)=(a_n,\ldots,a_1)$\\

When $a_n=b_n=1$ the base case gives the desired conclusion.  If $a_n=2$ we expand twice by recursion to obtain:
\begin{align*}
{\langle{e_{({b_n},\ldots,{b_1})}},{e'_{({a_n},\ldots,{a_1})}}\rangle}&={\langle{l_2(e_{({b_{n-1}},\ldots,{b_1})})},l_{2}(e'_{({a_{n-1}},\ldots,{a_1})})\rangle}\\
&-{\frac{\Delta_{0}}{\Delta_{1}}}{\cdot}{\langle{l_2(e_{({b_{n-1}},\ldots,{b_1})})},l_{1}(e'_{({a_{n-1}},\ldots,{a_1})})\rangle}\\
&=(q-{\frac{\Delta_{0}}{\Delta_{1}}}){\cdot}{\langle{e_{({b_{n-1}},\ldots,{b_1})}},{e'_{({a_{n-1}},\ldots,{a_1})}}\rangle}\\
&={\frac{\Delta_{2}}{\Delta_{1}}}{\cdot}{\frac{\Delta_{a_{n-1}}}{\Delta_{{a_{n-1}}-1}}}{\cdots}{\frac{\Delta_{a_{1}}}{\Delta_{{a_{1}}-1}}}\\
\end{align*}
which is the desired expression.
For $3{\leq}a_n$ we again expand twice, but an extra term results. 
\begin{align*}
&={\langle{l_{a_n=b_n}(e_{({b_{n-1}},\ldots,{b_1})})},l_{a_n=b_n}(e'_{({a_{n-1}},\ldots,{a_1})})\rangle}\\
&-{\frac{\Delta_{{a_n}-2}}{\Delta_{{a_n}-1}}}{\cdot}{\langle{l_{b_n=a_n}(e_{({b_{n-1}},\ldots,{b_1})})},l_{{a_n}-1}(e'_{({a_{n-1}},\ldots,{a_1})})\rangle}\\
&+{\frac{\Delta_{{a_n}-3}}{\Delta_{{a_n}-1}}}{\langle{e_{({b_{n}},\ldots,{b_1})}},e'_{({a_{n-2}},\ldots,{a_1})}\rangle}\\
\intertext{Setting $j={a_n}-2$ we have $j={a_n-2}={b_n}-2<{{b_n}-1}={{a_n}-1}{\leq}{a_{n-1}}$ we may apply Lemma~$\ref{sub}$ to obtain:}
&={\langle{e_{({b_{n}},\ldots,{b_1})}},e'_{({a_{n-2}},\ldots,{a_1})}\rangle}=0\\
\intertext{The remaining terms simplify to}
&={\langle{l_{a_n=b_n}(e_{({b_{n-1}},\ldots,{b_1})})},l_{a_n=b_n}(e'_{({a_{n-1}},\ldots,{a_1})})\rangle}\\
&-{\frac{\Delta_{{a_n}-2}}{\Delta_{{a_n}-1}}}{\cdot}{\langle{l_{b_n=a_n}(e_{({b_{n-1}},\ldots,{b_1})})},l_{{a_n}-1}(e'_{({a_{n-1}},\ldots,{a_1})})\rangle}\\
&=(q-{\frac{\Delta_{{a_n}-2}}{\Delta_{{a_n}-1}}}){\langle{e_{({b_{n-1}},\ldots,{b_1})}},e'_{({a_{n-1}},\ldots,{a_1})}\rangle}\\
&={\frac{\Delta_{{a_n}}}{\Delta_{{a_n}-1}}}{\langle{e_{({b_{n-1}},\ldots,{b_1})}},e'_{({a_{n-1}},\ldots,{a_1})}\rangle}\\
\intertext{and by induction on $n$ we obtain what we want}
&={\frac{\Delta_{{a_n}}}{\Delta_{{a_n}-1}}}{\cdot}{\frac{\Delta_{{a_{n-1}}}}{\Delta_{{a_{n-1}}-1}}}{\cdots}{\frac{\Delta_{{a_1}}}{\Delta_{{a_1}-1}}}\\
\end{align*}
\end{proof}
\end{section}

\begin{section}{The Markov form in the Orthogonal Basis}

We now prove our main result that the space of chord diagrams, $V_n$,
equipped with the Markov symmetric bilinear form, $\langle\ ,
\rangle_n$, has $\{e'_\alpha\}$ for an orthogonal basis.
\begin{theorem} The basis $\{e'_\alpha\}$ is orthogonal. Furthermore, for  $\alpha = (a_n,\ldots,a_1).$
\begin{itemize}
\item ${\langle{e'_{\alpha}},{e'_{\alpha}}}\rangle=\frac{\Delta_{a_n}}{\Delta_{{a_n}-1}}{\cdots}\frac{\Delta_{a_1}}{\Delta_{{a_1}-1}}$
\item If $\beta \ne \alpha$ then ${\langle{e'_{\beta}},{e'_{\alpha}}\rangle}=0$
\end{itemize}
\end{theorem}
\begin{proof}
By symmetry, we may assume
$\beta=(b_n,\ldots,b_1){\leq}\alpha=(a_n,\ldots,a_1)$. By the
corollary, we can write
${e'_{(b_n,\ldots,b_1)}}=\sum_{\alpha{\preccurlyeq}{\in}D_n}P_{\alpha}(q){\alpha}$
where the $P_{\alpha}(q)$ are rational polynomials which are 0 for
$\alpha{\npreceq}{(b_n,\ldots,b_1)}$.

Hence the computation of the Markov form yields:
${\langle{e'_{(b_n,\ldots,b_1)}},{e'_{(a_n,\ldots,a_1)}}\rangle}=\sum_{{\alpha{\preccurlyeq}(b_n,\ldots,b_1)}D_n}P_{\alpha}(q){\langle{\alpha},{e'_{(a_n,\ldots,a_1)}}\rangle}$.

In the case $(b_n,\ldots,b_1){<}(a_n,\ldots,a_1)$,  we have for
every $\alpha$ in the sum:
$\alpha{\preccurlyeq}(b_n,\ldots,b_1)\Rightarrow{\alpha}{\leq}(b_n,\ldots,b_1)\Rightarrow{\alpha}{<}{(a_n,\ldots,a_1)}$.
So by Theorem~$\ref{Fact Two}$, each bracket term is 0, so the sum
is 0.

When  $(b_n,\ldots,b_1){=}(a_n,\ldots,a_1)$ the same reasoning
applies, except that  there is a single term left, 
the one term where equality holds:
\[{P_{\alpha=(a_n,\ldots,a_1)}(q)}{\langle{e_{(a_n,\ldots,a_1)}},{e'_{(a_n,\ldots,a_1)}}\rangle}\]
But, as previously remarked, $P_{(a_n,\ldots,a_1)}(q)=1$ and by
$\ref{Fact Two}$ only the desired product
$\frac{\Delta_{a_n}}{\Delta_{{a_n}-1}}{\cdots}\frac{\Delta_{a_1}}{\Delta_{{a_1}-1}}$ remains.
\end{proof}


\end{section}\begin{section}{Relation with other orthogonal bases}

The following technique for constructing orthogonal bases for quantum modules, implicit in the literature, was shown to the authors by Pat Gilmer.  The method uses the calculus of trivalent graphs for computing three-manifold quantum invariants (see Kauffman-Lins\cite{KL94}).

The method begins with a fixed choice of a trivalent plane tree in the unit disk with $2n$ terminal (external) edges and $2n-3$ internal edges. In the method of Kauffman and Lins, the edges must be \textit{coloured} with (quantum) integers satisfying all the triangular inequalities (at each trivalent vertex the sum of the integers must be even and the sum of any two integers minus the third must be non-negative).  For this construction all the external edges must be coloured one.  The number of such colorings is the Catalan number, as expected.

Next,  we pair two such colorations of our chosen trivalent plane tree by gluing the univalent ends and computing the resulting quantum invariant in $\mathbb{Z}[q^\epsilon]$.  Using  two-spheres meeting the paired configuration in two points labelled one, we inductively compute the Markov pairing using the fusion rule (Lemma 7, p. 40\cite{KL94}). The components of the sum either all vanish (when the colorations are different) or 
one single summand remains, which is a  product of theta coefficients, easily rewritten as a quotient of Chebyschev polynomials. 

For $n = 3$, there are three non-isomorphic configurations  (removing the six terminal edges from a plane trivalent tree on 9 edges leaves a tree with three edges which must either be a linear chain or a trivalent \:\textbf{\large Y}\:. The \:\textbf{\large Y}\:  configuration completes to a unique plane trivalent tree  with 9 edges and the chain completes to two plane trivalent trees, differening in whether or not, the terminal edges connected to the interior vertices of the chain are on the same or opposite  \textit{sides}). The associated change of basis to the \textit{non-crossing} chord diagram basis (after making a choice of root for the chord diagram) and the corresponding orthogonal matrices are given below.

\subsection{Colored Trivalent Graphs}

First, for the \:\textbf{\large Y}\: configuration, we provide the change of basis matrix:
 \[
 \left(
   \begin{array}{lllll}
    1 & 0 & 0 & 0 & 0 \\
    -\frac{1}{q} & 0 & 1 & 0 & 0 \\
    -\frac{1}{q} & 0 & 0 & 0 & 1 \\
    -\frac{1}{q} & 1 & 0 & 0 & 0 \\
    \frac{2}{q^2} & -\frac{1}{q} & -\frac{1}{q} & 1 & -\frac{1}{q}
   \end{array}
   \right)
\]
and the Markov form in this orthogonal basis.
\[
\left(
                  \begin{array}{lllll}
                   q^3 & 0 & 0 & 0 & 0 \\
                   0 & (q-1) q (q+1) & 0 & 0 & 0 \\
                   0 & 0 & (q-1) q (q+1) & 0 & 0 \\
                   0 & 0 & 0 & (q-1) q (q+1) & 0 \\
                   0 & 0 & 0 & 0 & \frac{(q-1) (q+1) \left(q^2-2\right)}{q}
                  \end{array}
                  \right)
\]

Next the same information is provided for the \textit{same side} configuration:

\[
  \left(
   \begin{array}{lllll}
    1 & 0 & 0 & 0 & 0 \\
    -\frac{1}{q} & 1 & 0 & 0 & 0 \\
    -\frac{1}{q} & 0 & 1 & 0 & 0 \\
    \frac{1}{q^2} & -\frac{1}{q} & -\frac{1}{q} & 1 & 0 \\
    -\frac{q}{q^2-1} & \frac{1}{q^2-1} & \frac{1}{q^2-1} & -\frac{q}{q^2-1} &
      1
   \end{array}
   \right)
\]
\[
\left(
                  \begin{array}{lllll}
                   q^3 & 0 & 0 & 0 & 0 \\
                   0 & (q-1) q (q+1) & 0 & 0 & 0 \\
                   0 & 0 & (q-1) q (q+1) & 0 & 0 \\
                   0 & 0 & 0 & \frac{(q-1)^2 (q+1)^2}{q} & 0 \\
                   0 & 0 & 0 & 0 & q \left(q^2-2\right)
                  \end{array}
                  \right)
\]

and, finally  for  the opposite side configuration:

\[
\left(
\begin{array}{lllll}
0 & 0 & 1 & 0 & 0\cr     
1 & 0 & - \frac{1}{q}   & 0 & 0\cr     
0 & 0 & -\frac{1}{q}  & 1 & 0\cr
- \frac{1}{q}  & 1 & q^{-2} & - \frac{1}{q}   & 0\cr 
- \frac{q}{-1 + q^2}   & \frac{1}{-1 + q^2} & \frac{1}{-1 + q^2} & \frac{q}{-1 + q^2}   & 1
                  \end{array}
\right)
\]

\[
\left(
                  \begin{array}{lllll}
                   q^3 & 0 & 0 & 0 & 0 \\
                   0 & (q-1) q (q+1) & 0 & 0 & 0 \\
                   0 & 0 & (q-1) q (q+1) & 0 & 0 \\
                   0 & 0 & 0 & \frac{(q-1)^2 (q+1)^2}{q} & 0 \\
                   0 & 0 & 0 & 0 & -2  q + q^3
                  \end{array}
                  \right)
\]

\subsection{Computation for the poset recursion construction}
Here is the change of basis matrix for the poset-recursion.

\[
   \left(
   \begin{array}{lllll}
    1 & 0 & 0 & 0 & 0 \\
    -\frac{1}{q} & 1 & 0 & 0 & 0 \\
    -\frac{1}{q} & 0 & 1 & 0 & 0 \\
    \frac{1}{q^2} & -\frac{1}{q} & -\frac{1}{q} & 1 & 0 \\
    -\frac{q}{q^2-1} & \frac{1}{q^2-1} & \frac{1}{q^2-1} & -\frac{q}{q^2-1} &
      1
   \end{array}
   \right)
\]
and the matrix for the Markov form.

\[
\left(
                  \begin{array}{lllll}
                   q^3 & 0 & 0 & 0 & 0 \\
                   0 & (q-1) q (q+1) & 0 & 0 & 0 \\
                   0 & 0 & (q-1) q (q+1) & 0 & 0 \\
                   0 & 0 & 0 & \frac{(q-1)^2 (q+1)^2}{q} & 0 \\
                   0 & 0 & 0 & 0 & q \left(q^2-2\right)
                  \end{array}
                  \right)
\]
The \textit{same-side} configuration gives the same diagonal form and triangular change of basis as the recursively computed diagonalized matrix.
Although the orthogonal basis is the same,   the poset-recursion formula provides an explanation for the location of the non-zero coefficients in the triangular change of basis:  the non-zero coefficients can only occur in positions corresponding to those restricted sequences smaller than the given one, that is, in the down-set of the given element in the coordinate-wise partial order.


\end{section}\begin{section}{Further Development}

In a subsequent paper, we will develop further properties of this poset from the perspective of \textit{non-crossing fixed point free involutions} in a permutation group of a set of even cardinality.  In addition, this new viewpoint will allow us to extend this theory to annular braids, that is, Artin braids of type $\mathbb{B}_n$.

\end{section}

\end{document}